\documentclass[a4paper]{amsart}
\usepackage{amssymb,latexsym,epsfig}
\usepackage[numbers,sort]{natbib}
\usepackage{amsmath}
\usepackage{amsfonts}
\usepackage{graphicx}

\def\restr#1_#2{\setbox1=\hbox{$#1$}\setbox2=\hbox{$_{#2}$}\dimen1=\ht1%
\dimen2=\dp2\box1\vrule width 0.4pt height \dimen1 depth \dimen2\mskip1mu \box2}%
\let\Schrstr=/%
\def \fff{\hbox to 0pt{\hss}}%

\def\drauf#1#2{{\setbox1=\hbox{#1}\setbox2=\hbox{#2}\dimen1=\wd1\fff%
\ifdim\dimen1<\wd2\dimen1=\wd2\else\fi\hbox to \dimen1{\hfill#1\hfill}%
\hskip-\dimen1\hbox to\dimen1{\hfill#2\hfill}}}%

\newdimen\epsfxsize
\newbox\epsfbox
\setbox\epsfbox=\hbox{~~~}%

\theoremstyle{plain}

\newtheorem{theo}{Theorem}[section]

\newtheorem{prop}[theo]{Proposition}

\newtheorem{defn}[theo]{Definition}
\newtheorem{rmrk}[theo]{Remark}

\makeatletter

\def\a{\alpha}

\def\g{\gamma}
\def\d{\delta}

\def\r{\rho}

\def\Komma{$\discretionary{\hbox{~},}{}{\hbox{~~},\hbox{~~}}$}%
\def\invlim{\vtop{\offinterlineskip\hbox{lim}\hbox{$\longleftarrow$}}}%
\def\calU{{\mathcal U}}%
\def\calV{{\mathcal V}}%
\def\RR{\mathbb{R}}

\def\ZZ{\mathbb{Z}}
\def\UU{\mathbb{U}}

\let\noi=\noindent
\begin{document}

\bigskip

\title[On semilocally simply connected spaces]{On semilocally simply connected spaces}

 \author{Hanspeter Fischer}
 \address{Department of Mathematical Sciences,
          Ball State University, Muncie, IN 47306, United States}
\email{fischer@math.bsu.edu}

\author{Du\v san Repov\v s}
\address{Faculty of Mathematics and Physics,
and Faculty of Education,
University of \newline Ljubljana,
Jadranska 19,
Ljubljana 1000,
Slovenia}
\email{dusan.repovs@guest.arnes.si}

\author{\v Ziga Virk}
\address{Faculty of Mathematics and Physics,
University of \newline Ljubljana,
Jadranska 19,
Ljubljana 1000,
Slovenia}
\email{zigavirk@gmail.com}

\author{Andreas Zastrow}
\address{Institute of Mathematics, Gdansk University, ul. Wita Stwosza 57,
80-952 Gda\'nsk, Poland}
\email{zastrow@mat.ug.edu.pl}

\date{\today}

\keywords{Semilocal simple connectivity, generalized universal covering space, fundamental group, Spanier group, homotopically Hausdorff spaces, shape injectivity, inverse limit}
\subjclass[2010] {Primary: 54D05;  Secondary: 55Q05, 54G20, 54G15,  57M10, 55Q07}

\begin{abstract}
The purpose of this paper is:
(i) to construct a space which is semilocally simply connected
    in the sense of Spanier even though
    its Spanier group is non-trivial;
(ii) to propose a modification of the notion of a Spanier group
so that via the modified Spanier group
semilocal simple connectivity can be
characterized; and
(iii) to point out that with just a slightly modified definition of semilocal simple connectivity
which is sometimes also used in literature,
the classical Spanier group gives the correct
characterization within the general class of path-connected topological spaces.\par
While the condition ``semilocally simply connected"
plays a crucial role in classical covering theory,
in generalized covering theory one needs to consider the condition ``homotopically Hausdorff\/\/" instead.
The paper also discusses which implications hold between all
of the abovementioned conditions and, via the modified Spanier groups, it also
unveils the weakest so far known
algebraic characterization for the existence
of generalized covering spaces
as introduced by Fischer-Zastrow.
For most of the implications,
the paper also proves the non-reversibility
by providing the corresponding examples.
Some of them rely on spaces that are newly
constructed in this paper.
\end{abstract}

\maketitle

\section{Introduction}

This paper was motivated by an observation during the research of
\cite{FZ}, namely that E.H. Spanier, when writing his celebrated
book on algebraic topology \cite{Sp},
apparently made an oversight in the statement which
immediately precedes
Corollary 2.5.14. That statement, in which he characterizes semilocal
simple connectivity in terms of vanishing of a certain group $\pi({\mathcal U},x_0)$ for at least one open covering $\mathcal U$ of the space,
turns out to be correct only if one
additionally assumes local path-connectedness. Of course, one may wonder
if this assumption was perhaps not implicitly made. However, in view of the
author's great attention to details in this book in general, our
ultimate conclusion was that Spanier would have
mentioned this additional assumption had he been aware of the phenomena and
examples that we shall expose below.

Roughly speaking, a problem occurs for spaces which are not locally
path-connected. This is because the fundamental group uses
base points, whereas the subgroup of the fundamental  group which
Spanier
associated with a covering (and which we henceforth call the {\it Spanier group})
does not use base points in a similar way. Of course, this group is
defined as a subgroup of the fundamental group and therefore it depends on
the same  base point as the fundamental group, but the concept of
Spanier groups refrains from using base points for each set of the covering
separately. Therefore Spanier's characterization of semilocal
simple connectivity (cf.\ immediately before Corollary 2.5.14)
matches his definition (cf.\ immediately before Theorem 2.4.10) only in the
locally path-connected case.

The main purpose of our paper is to:
\begin{itemize}
\item   confirm this assertion by constructing a space $Y$ which is
        semilocally simply connected in the sense of Spanier, but its
        Spanier group is non-trivial (cf.\
        Proposition \ref{Y});
\item   confirm the ``if\/\/"-part of the claim preceding
        Corollary 2.5.14, for all spaces and the ``only-if\/\/" part for locally path-connected
        spaces
        (Theorem \ref{properties}(3)--(4));
\item   propose a modification of Spanier groups so that the corresponding claim will
        be correct for all spaces
        (Theorem \ref{properties}(2)) and
\item   propose a modification of the definition of semilocal simple connectivity
        and prove that with this modified definition and with the original
        definition of Spanier groups the claim preceding Corollary 2.5.14 is also correct
        (Theorem \ref{properties}(1)).
\end{itemize}
\

Accordingly, there will be two concepts of semilocal simple connectivity
and two versions of Spanier groups --- one which depends on base points and one which does not.
In order to avoid  ambiguity in the terminology of the present paper, we will
from now on speak of these concepts using the attributes ``based" and
``unbased'' (cf.\ Definitions \ref{BasedSLSC}--\ref{UnbasedSpanierGroup} and \ref{BasedSpanierGroup}).

While semilocal simple connectivity is a crucial condition in classical
covering space theory, the generalized covering space theory  treated
in \cite{FZ} mainly considered the weaker condition called
``homotopically Hausdorff\/\/".
The paper \cite{CMRZZ} also studied two versions of this condition,
one which depends on base points and one which does not.
 However, in our case we will
adopt the notation from \cite{CMRZZ} and use the attributes ``weak" and ``strong" (cf.\ Definition \ref{homHaus}),
respectively. When base points are treated correctly, the conditions
of semilocal simple connectivity can be equivalently described by properties
of Spanier groups (Theorem \ref{properties}(1)--(2)). However, we only know sufficient conditions on Spanier groups
which imply homotopic Hausdorffness
(Section 6, (11)--(12)). In this context it should
also be pointed out that even the weakest of these conditions for Spanier groups
implies a condition that we will call ``homotopically path-Hausdorff\/\/".
This condition was not mentioned in \cite{CMRZZ}
or \cite{FZ}, but is similar to a condition that appeared under a different name
in \cite{Z-old}, a preprint preceding \cite{FZ}.
Therefore the concept of based Spanier groups, introduced in this paper,
also apparently yields the weakest currently known algebraic sufficient
condition for the existence of generalized universal covering spaces
(Theorem \ref{covering}).

We also wish to point out that we are not aware of the above
described incorrectness in Spanier's book leading to a false theorem therein.
All places that we found, where the crucial remark preceding
Corollary 2.5.14 has been applied, were statements where  local
path-connectedness of the underlying topological space has been an
assumption, and in this framework the crucial statement is correct.

Our paper will in Section~6 also briefly discuss known implications among all
the mathematical concepts mentioned so far.
All topological spaces in this article are assumed to be path-connected.
Since we need to consider the fundamental group, it does not seem necessary
to consider more general spaces.

\section{Definitions and Terminology}

\begin{defn}\label{BasedSLSC}
We call a topological space $X$
{\rm (based) semilocally simply connected}
if for every point $x\in X$
there exists a neighbourhood $U$ of $x$
such that the inclusion-induced homomorphism
$\pi_1(U,x)\to \pi_{1}(X,x)$
is trivial.
\end{defn}

The majority of topology books discussing covering spaces
seems to prefer this definition (see e.g.\ p.\ 63 of \cite{Ha}, p.\ 393 of \cite{Mu}, p.\ 174 of \cite{M}
and p.\ 78 of \cite{Sp}) over the following; we will only
use the attribute ``based'' in connection
with semilocal simple connectivity where it is needed to distinguish
this definition from the following.

\begin{defn}\label{UnbasedSLSC}
We call a topological space $X$
{\rm unbased semilocally simply connected}
if for every point $x\in X$
there exists a neighbourhood $U$ of $ x$
such that every loop in $U$ is null-homotopic in $X$.
\end{defn}

The latter definition is used, for example, in \cite[p.\ 187]{F},
and, named slightly differently, in \cite[Definition 6.6.8, p.\ 255]{HiWy}.

\begin{defn}\label{UnbasedSpanierGroup}
Let $X$ be a space, $x_{0} \in X$ a base point,
and ${\mathcal U} = \{U_i \mid i \in I \}$
an arbitrary open covering of $X$.
Then  we define $\pi({\mathcal U}, x_{0})$
to be the subgroup of
$\pi_{1}(X,x_{0})$
which contains all homotopy classes
having representatives of the following type:
$$\prod^n_{j=1} u_j~v_j~u_j^{-1}, \leqno (\ref{UnbasedSpanierGroup})$$
where $u_j$ are arbitrary paths (starting at the base point $x_0$)
and each $v_j $ is a
loop
inside one of the neighbourhoods
$U_i \in {\mathcal U}$.
We call this group the {\rm (unbased) Spanier group with respect to ${\mathcal U}$}.
\end{defn}

This definition matches the definition from \cite[Chapter 2, Section 5 between items 7 and 8]{Sp}.
In \cite{Sp} this notation was already used, but names have not yet
been given to these groups.

We choose the name {\it Spanier
groups}, since all traces in the literature that we are aware of seem to go back to this
appearance in Spanier's book.

In the introduction we announced  a concept that introduces base points
to the sets of open coverings. {Consequently, we will instead
of open sets $U$ also consider {\it``pointed open sets''}, i.e.\
pairs $(U,x)$, where $x\in U$ and $U$ is open.

\begin{defn}\label{nbd-pair} Let $X$ be a space.
\begin{enumerate}
\item An {\rm open covering of $X$ by pointed sets} is
      a family of pointed open sets ${\mathcal V} = \{ (U_i,x_i) \mid i \in I\}$,
      where
$$                             \{x_i \mid i \in I\}   = X.\leqno \hbox{\rm\hskip2.9pc(\ref{nbd-pair})}$$
\item {\rm Refinements} between coverings
by pointed sets are defined as follows:
                      ${\mathcal U}' = \{ (U'_i,x'_i )  \mid i \in I\}$ refines
                     $ {\mathcal U}  = \{ (U_j, x_j)   \mid j \in J \}$ ,
               if $\forall i \; \exists j$ such that $U'_i \subset U_j $ and $x_i = x_j$.
\end{enumerate}
\end{defn}

Let ${\mathcal U} = \{ U_i \mid i \in I\}$
be a covering of $X$ by open sets. Observe that due to expression (\ref{nbd-pair}) demanding that each
point of $X$ occurs at least once as base point of one of the covering sets,
it will in general not suffice to  choose a base point for each
of the $U_i$ in order to turn it into an open covering ${\mathcal V}$ of $X$ by pointed sets.
Instead, the following procedure is apparently in general necessary:

\begin{itemize}
\item for each $U_i \in {\mathcal U}$ take  $|U_i|$ copies into ${\mathcal V}$; and
\item define each of those copies as $(U_i,P)$, i.e.\ use the same set $U_i$
      as first entry, and let the second entry run over all points $P \in U_i$.
\end{itemize}

When constructed with this procedure, coverings by neighbourhood pairs offer in principle the same options
for refinements as coverings by open sets.
Vice versa note, that this procedure will usually generate such coverings
by pointed sets, where a lot of $x \in X$ occur as base points
for different sets $U_i$}.

\begin{defn} \label{BasedSpanierGroup}{%
Let $X$ be a space, $x_{0}\in X$, and ${\mathcal V} = \{(U_i,x_i) \mid i \in I \}$
be a covering of $X$ by open neighbourhood pairs. Then
we let $\pi^*({\mathcal V}, x_{0})$ be the subgroup of $\pi_1(X,x_0)$ which contains all
homotopy classes
having representatives of the following type:
$$\prod^n_{j=1} u_j~v_j~u_j^{-1},\leqno(\ref{BasedSpanierGroup}) $$
where the $u_j$ are arbitrary paths that run from $x_0$ to some point $x_i$
and each $v_j$ then must be a closed path inside the corresponding $U_i$}.
We will call this group the {\rm based Spanier group with respect to ${\mathcal V}$}.
\end{defn}

\begin{rmrk}
\label{remark0}
Note that for based and unbased Spanier groups the following holds:
Let ${\mathcal U} , {\mathcal V}$ be open coverings, and let ${\mathcal U}$ be a
refinement of ${\mathcal V}$.
Then $\pi({\mathcal U}) \subset \pi({\mathcal V})$.
Analogously $\pi^*({\mathcal U}) \subset \pi^*({\mathcal V})$
holds, assuming that ${\mathcal U}$ and ${\mathcal V}$ are now open coverings by pointed sets.
Due to these inclusion relations, there exist inverse limits of these Spanier
groups, defined via the directed system of all coverings with respect to
refinement. We will call them {\rm the (unbased) Spanier group and the based Spanier group of
the space $X$}
and denote them by
{\rm
$$\invlim(\pi({\mathcal U}))\hbox{\it~~ and~~} \invlim(\pi^*({\mathcal V})),$$ }respectively,
observing that these inverse limits are realized by intersections:
{\rm$$\bigcap_{\text{coverings }\mathcal U} \pi({\mathcal U}) =
\invlim(\pi({\mathcal U}))\hbox{\it~~ and~~}
\bigcap_{\text {coverings }\mathcal V \text{ by pointed sets}} \pi^*({\mathcal V})
= \invlim(\pi^*({\mathcal V})).$$
}%
Since for locally path-connected
spaces any covering refines to a covering by path-connected sets,
the based and unbased Spanier groups
coincide in that case.

\end{rmrk}

\begin{rmrk}
\label{remark}
Recall that all topological spaces in this paper are assumed to be
path-connected. Despite  discussing the properties ``based"
and ``unbased", by definition the fundamental group and (since they are defined
as subgroups) also all Spanier groups formally depend on a base point. However,
the standard argument that the isomorphism type and many other essential properties
of the fundamental group of a path-connected space do indeed not depend
on the base point, naturally extends to Spanier groups. We will therefore
omit the base point in our notation, whenever appropriate.
\end{rmrk}

The following two theorems are our main results, and Sections~4~ and ~5,
respectively, are devoted to prove them.

\begin{theo}\label{properties}\hspace{0pt}
\begin{enumerate}
\item Let $X$ be an arbitrary space, $x_{0} \in X$. Then $X$ is unbased semilocally simply connected,
    if and only if it has an open covering ${\mathcal U}$ such that $\pi({\mathcal U},x_{0})$ is trivial.
\item Let $X$ be an arbitrary space, $x_{0} \in X$. {Then $X$ is  semilocally simply connected,
    if and only if it has an open covering ${\mathcal V}$ by pointed sets}
    such that $\pi^*({\mathcal V},x_{0})$ is trivial.
\item Let $X$ be an arbitrary space. Then the two equivalent properties from (1)
    imply those from (2).
\item Let $X$ be a locally path-connected  space. Then the two equivalent properties from (2)
    imply those from (1).
\item For topological spaces that are not locally path-connected, the properties
    from (2) need not imply those from (1). The space $Y$ that will be constructed
    in Section ~3, satisfies (2) but not (1).
\end{enumerate}
\end{theo}

\begin{theo}\label{covering}
Let $X$ be a topological space whose based Spanier group
$\invlim\pi^*({\mathcal U})$ is trivial. Then $X$ admits a generalized universal covering space
in the sense of \cite{FZ}.
\end{theo}

With items (1) and (2) of the next definition we follow
\cite[p.\ 1091]{CMRZZ}:

\begin{defn} \label{homHaus}\hspace{0pt}
\begin{enumerate}
\item  A space $X$ is called (weakly) {\rm homotopically Hausdorff} if for every $x_0 \in X$ and for every  non-trivial $\alpha \in \pi_1(X,x_0)$ there exists a neighbourhood $U$ of $x_0$ such that no loop in $U$ is homotopic (in $X$) to $\alpha$ rel.\ $x_0$. An equivalent condition, using the terminology of \cite[Definition 1]{ZV}, would be the absence of non-trivial small loops.
\item  A space $X$ is called {\rm strongly homotopically Hausdorff}, if for every
   $x_0 \in X$ and for every essential closed curve $\gamma$ in $ X$ there is a neighbourhood of $x_0$ that contains no closed curve freely homotopic (in $X$) to $\gamma$.
\item A space $X$ is called {\rm homotopically path-Hausdorff},
provided that it satisfies the following property with respect to any
two paths  $w,v : [0,1] \to X$ with $w(0) = v(0)$ and $w(1) = v(1)$: If
$w$ and $v$ are not homotopic relative to the endpoints, then there exist
$0=t_0<t_1<\ldots<t_k=1$ and open sets $U_1,U_2,\ldots, U_k$ with
$w([t_{i-1},t_i]) \subseteq U_i$ for $1 \leq i \leq k$ such that for any
$w_i :  [t_{i-1},t_i] \to U_i$ with $w_i(t_{i-1})\mathord{ =} w(t_{i-1})$ and
$w_i(t_{i})\mathord{ = }w(t_{i})$, the concatenation $w_1* w_2*\ldots* w_k$ is not
homotopic to $v$ relative to the endpoints.

\end{enumerate}
\end{defn}
 The terms weakly and strongly homotopically Hausdorff
already appeared, but under different names, in \cite[1.1]{Z-old},
while the third property that was considered there, although similar to our
Definition \ref{homHaus}(3), is not equivalent to it.
While `strongly
homotopically Hausdorff\/' appeared for the first time in peer-reviewed
literature in \cite{CMRZZ}, `homotopically Hausdorff\/', had been used before in \cite[Definition 5.2]{CC}.

\begin{rmrk}\label{homHausRelation}
{\rm Note that all strongly homotopically Hausdorff spaces and all  homotopically path-Hausdorff spaces are homotopically Hausdorff.}
If we
apply the condition of homotopical path-Hausdorffness to a constant path at $x$ we obtain the condition for weak homotopical Hausdorffness at $x$;
thus the second statement. The first statement directly follows
from the definitions.

\end{rmrk}

\section{Examples}

\begin{figure}
\begin{center}
\includegraphics[height=8cm]{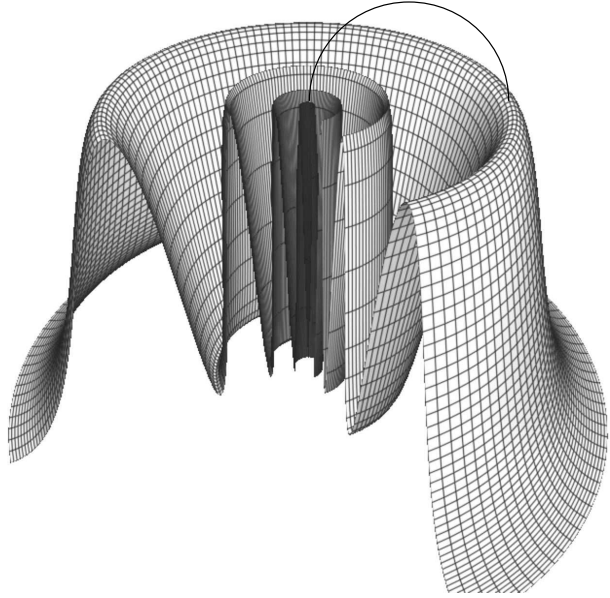}
\end{center}
\caption{ The space $Y$.}
\label{inn}
\end{figure}

\begin{figure}
\includegraphics{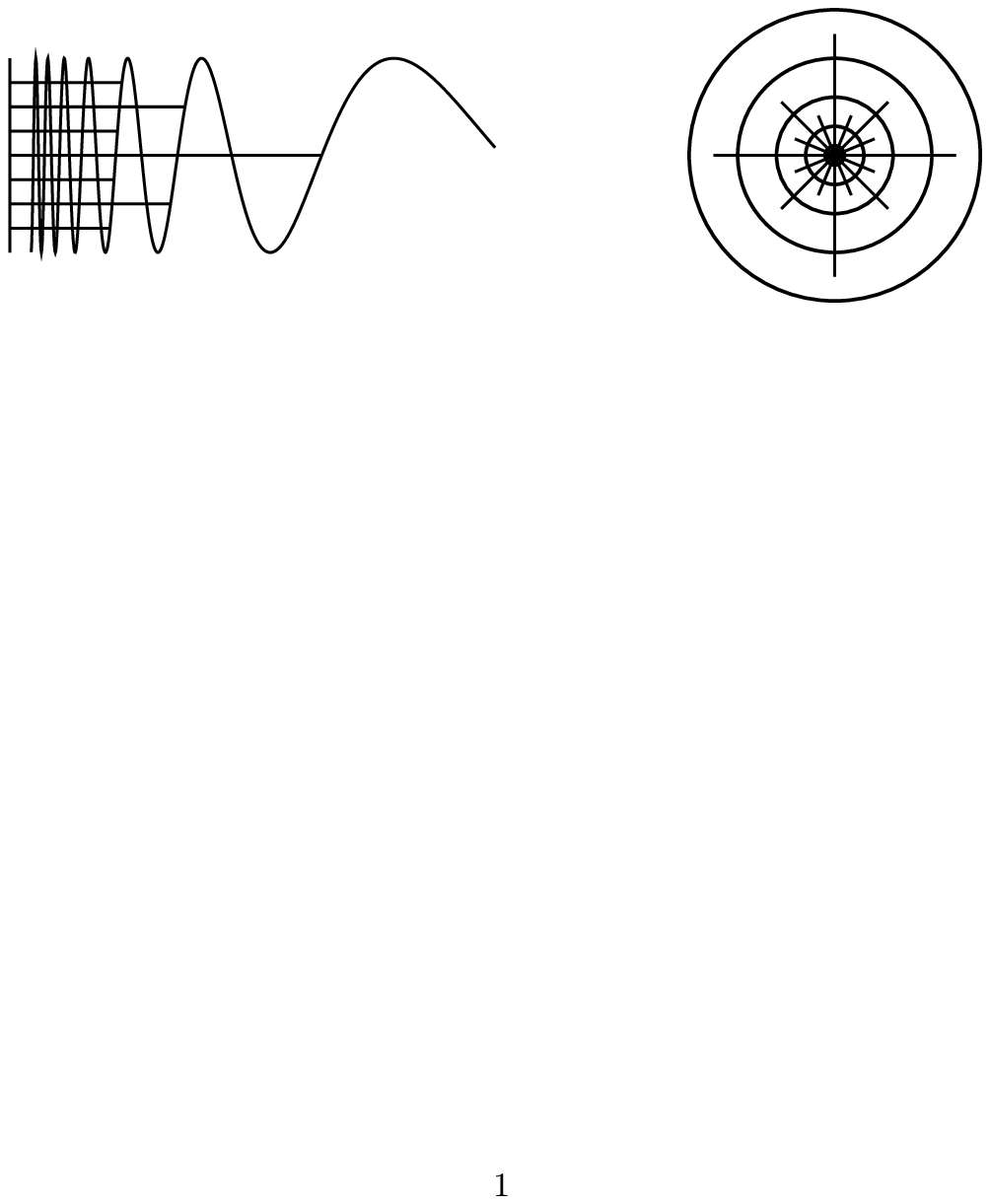}
\caption{The radial projections and the top view of the tunnels of $Y'$.}
\label{inn1}
\end{figure}

\begin{figure}
\begin{center}
\vbox to 8cm{\vss\epsfig{file=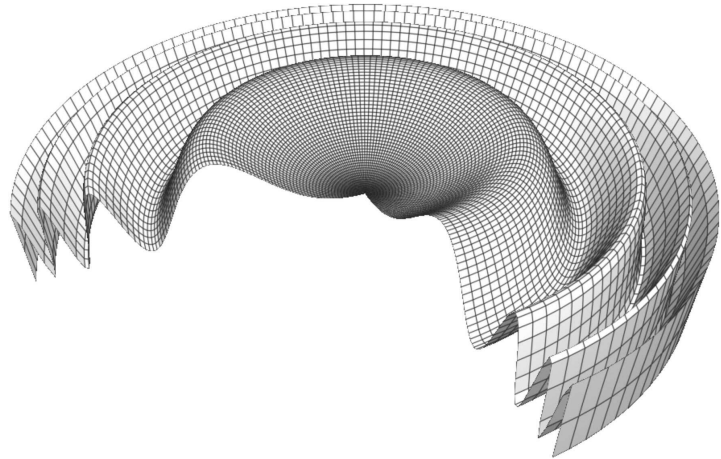,height=2.6in}}
\end{center}
\caption{ The ``surface'' portion of the spaces $Z$ and $Z'$.}
\label{out}
\end{figure}

\begin{figure}
\begin{center}
\epsfig{file=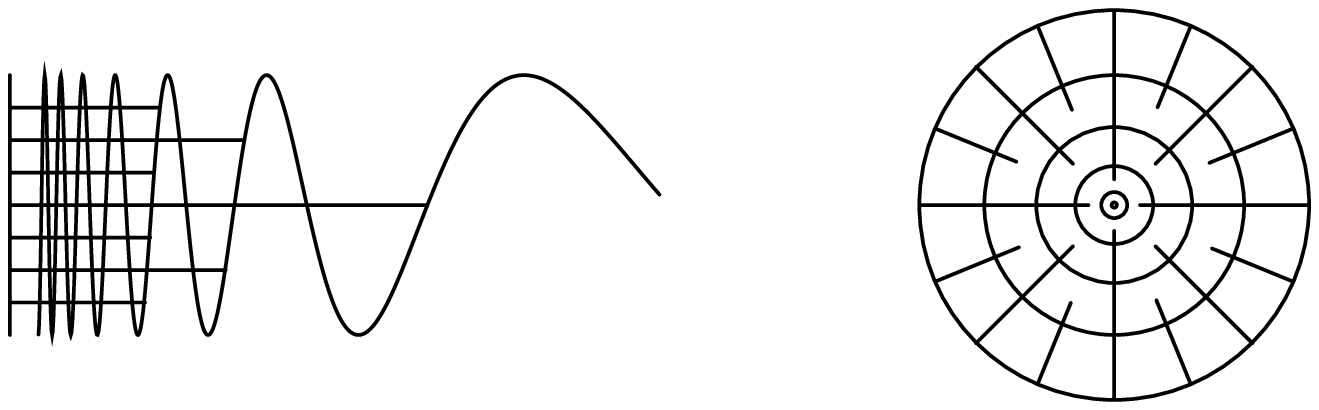,width=4in}
\end{center}
\caption{The radial projections and the top view of the tunnels of $Z'$.}
\label{out1}
\end{figure}

\begin{description}
\item[The space $Y'$] One of our examples $Y'\subset \RR^3$ will be  precisely the space
    called $A$ in \cite{CMRZZ}, and was defined
    at the beginning of Section 3 therein.  It consists of a rotated topologists' sine curve (as suggested by Figure \ref{inn}),
    the ``central axis'', where this surface tends to, and a system
    of horizontal arcs is attached to them so that they
    become dense (only) near the central axis.
    Figure \ref{inn1} shows schematically the top-view and the side-view of defining
    such arcs, which we shall call ``tunnels''. Their presence is indicated by the prime of $Y'$.

\item[The space $Y$] Another important example will be our space $Y\subset \RR^3$ (see Figure \ref{inn}). It consists of the same surface portion as $Y'$, the central
    axis to which this surface portion tends, but instead of defining a system
    of arcs for connecting central axis and surface portion,  we just connect them by a single arc $C$. This arc $C$ can be easily embedded into $\mathbb R^3$, so as not to intersect the surface portion or the central axis at any other points than its endpoints.

\item[The space $Z'$] Our third example will be called $Z'\subset \RR^3$. It  is precisely the space called $B$ in \cite{CMRZZ} and defined immediately before Theorem 3.4 therein. It consists of a rotated topologists' sine curve (as shown in Figure \ref{out}), the ``outer cylinder'' at radius $1$, where this surfaces tends to, and a system of horizontal arcs attached to them so that they  become dense (only) near the outer cylinder. Figure \ref{out1} shows schematically the top-view and the side-view of how to define such arcs. We will call these arcs ``tunnels'' as well.

\item[The space $Z$] Analogously, we will also need a space which has the same outer cylinder and surface portion as $Z'$, but where outer cylinder and surface portion are just connected by a single arc $C$ (similar as for $Y$), but which  has no tunnel-system. We  will call this space $Z$.
\end{description}


\begin{figure}
\begin{center}
\epsfig{file=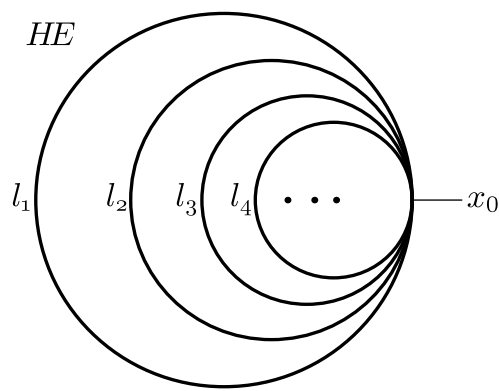,width=2.5in}
\end{center}
\caption{The Hawaiian Earring}
\label{HEfig}
\end{figure}

\noi Apart from  these spaces, we will also need the {\it Hawaiian Earring,}\/
that we will denote by $H\!\!E$.
This is a more well-known space, it is a countable union of
circles in the plane, as pictured in Figure \ref{HEfig}.

\begin{prop}\label{Y}
The space $Y$ has the following properties:
\begin{enumerate}
  \item Its Spanier group is non-trivial;
  \item its based Spanier group is trivial;
  \item it is semilocally simply connected;
  \item it is not unbased semilocally simply connected;
  \item it is homotopically Hausdorff;
  \item it is not strongly homotopically Hausdorff;
  \item it is homotopically path-Hausdorff.
\end{enumerate}
\end{prop}

\begin{prop}\label{Y'}
The space $Y'$ has the following properties:
\begin{enumerate}
  \item Its Spanier group is non-trivial;
  \item its based Spanier group is non-trivial;
  \item it is not semilocally simply connected;
  \item it is not unbased semilocally simply connected;
  \item it is homotopically Hausdorff;
  \item it is not strongly homotopically Hausdorff;
  \item it is homotopically path-Hausdorff.
\end{enumerate}
\end{prop}

\begin{prop}\label{Z}
The space $Z$ has the following properties:
\begin{enumerate}
  \item Its Spanier group is trivial;
  \item its based Spanier group is trivial;
  \item it is semilocally simply connected;
  \item it is unbased semilocally simply connected;
  \item it is homotopically Hausdorff;
  \item it is strongly homotopically Hausdorff;
  \item it is homotopically path-Hausdorff.
\end{enumerate}
\end{prop}

\begin{prop}\label{Z'}
The space $Z'$ has the following properties:
\begin{enumerate}
  \item Its Spanier group is non-trivial;
  \item its based Spanier group is non-trivial;
  \item it is not semilocally simply connected;
  \item it is not unbased semilocally simply connected;
  \item it is homotopically Hausdorff;
  \item it is strongly homotopically Hausdorff;
  \item it is not homotopically path-Hausdorff.
\end{enumerate}
\end{prop}

{\bf Convention:} Given any path $\a$ defined on  $[0,1]$, we define a path $\a^{-1}$ by $\a^{-1}(t):=\a(1-t)$.

{\bf Proof of Proposition \ref{Y}. }
(1) Fix a point $x_0$ on the surface portion of $Y$. Let $\r_r$ denote a simple path on the surface starting at $x_0$, contained in the plane determined by $x_0$ and the central axis, with endpoint at distance $r$ from the central axis.  Note that a simple loop $\a_r$ of constant radius $r>0$ on the surface is not freely homotopically trivial by \cite[Lemma 3.1]{CMRZZ}. Any neighbourhood of a point of the central axis contains such a loop. For every $1>r>0$ loops $\r_r~ \a_r ~ \r_r^{-1}$ (with $\a_r$ appropriately based) are homotopic to each other and non-trivial.  Given a cover of $Y$, all such loops are contained in the Spanier group of such a cover for sufficiently small $r$. Hence the Spanier group is non-trivial.

(2) Note that every point  $x\in Y$ has an arbitrarily small neighbourhood whose path component containing $x$ is contractible. Given a point on the surface at the radius $r$, an open ball of radius at most $r/2$ suffices. The interior of the connecting  arc $C$ induces such neighbourhoods for points of itself. Any open set not containing the entire arc $C$ suffices for the central axis, e.g.\ an open ball of radius at most $c/2$ where $c$ is the length of arc $C$.  The claim  follows by definition as the loops $v_j$ of Definition \ref{BasedSpanierGroup} are contractible. This also proves (3).

(4) The loops $\a_r$ from the proof of (1) are non-trivial and arbitrarily close to the central axis.

(5) Implied by (7) and Remark \ref{homHausRelation}.

(6) Follows from the proof of  (1) using the loops $\r_r$.

(7) We  use the notation of  Definition \ref{homHaus} and the fact that every point  $x\in Y$ has an arbitrarily small neighbourhood whose path component containing $x$ is contractible. Given any open cover $U_1, \ldots,U_k$ of $w([0,1])$ by such neighbourhoods, the only possible homotopy class for a path $w_1*\ldots*
w_k$ is that of $w$, constructing the homotopy using the contractibility of the sets $U_i$.
Thus the product of the $w_i$ will not be homotopic to $v$.
\hbox{~~~~}\hfill $\blacksquare$
\bigskip

{\bf Proof of Proposition \ref{Y'}. }
(1) The proof of  \ref{Y}(1) suffices.   Statements (1) and (2) are equivalent as $Y'$ is locally path-connected. Statements (3) and (4) are equivalent for the same reason and follow from (1).

(5) Follows from (7).

(6) Similarly as in proof of \ref{Y}  loops $\r_r$ provide an obstruction to strong homotopic Hausdorffness.

(7)
We will sketch a proof in \ref{hard}. Since the argument is
quite lengthy, a complete proof will appear elsewhere.
\hfill $\blacksquare$
\bigskip

{\bf Proof of Proposition  \ref{Z}. }
Note that every point  $x\in Z$ has an arbitrarily small neighbourhood whose every path component is contractible. For any point on the surface at radius $r$ an open ball of radius at most $(1-r)/2$ suffices. The interior of the connecting  arc $C$ induces such neighbourhoods for points of itself. Any open set not containing the entire arc $C$ suffices for the outer cylinder, e.g.\ an open ball of radius at most $c/2$ where $c$ is the length of the arc $C$.  Such cover proves (1)--(7). \hbox to\hsize{\hfill$\blacksquare$}
\bigskip

{\bf Proof of Proposition \ref{Z'}. }
(1) and (2) are equivalent by Theorem \ref{properties} and imply (3) and (4).

(1) For every positive $r<1$ let $\a_r$ denote a positively oriented simple loop  of constant radius  on the surface. All such loops are homotopically trivial. Let $\a_1$  be a clockwise oriented simple closed curve, defined as an intersection of the outer cylinder with the plane at height zero. Note that $\a_1$ is not  homotopically trivial by
\cite[Proof of Lemma 3.1]{CMRZZ}. We will prove that the loop $\a_1$ based at $x_0$ on the outer cylinder is contained in the Spanier group of $Z'$.

\begin{figure}
\includegraphics{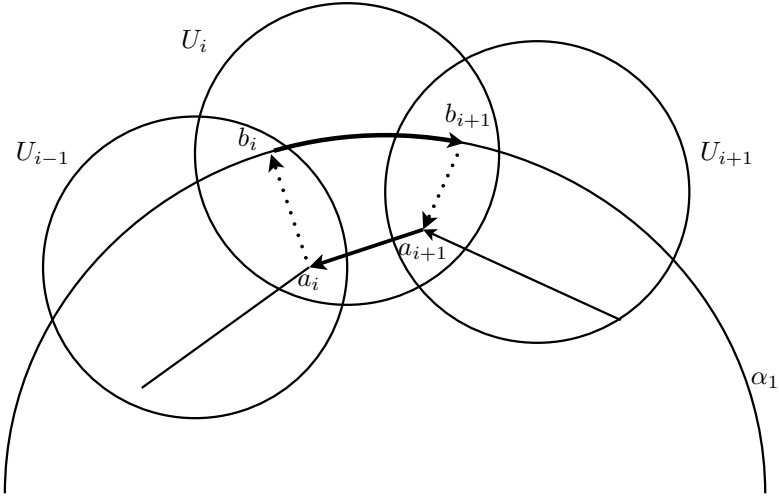}
\caption{Quadrilateral loop $Q_i$. The pointed lines mark path through the tunnel. }
\label{spaniergroupZ}
\end{figure}

Given any open cover $\calU$ of $\a_1([0,1])$ choose a finite refinement by balls $U_0, \ldots, U_k$ so that $U_i\cap U_j \neq \emptyset$ for $i,j\in \ZZ_{k+1}$ iff $|i-j|\leq 1$. Choose $R<1$ big enough so that $\a_R([0,1])$ is covered by $U_0, \ldots, U_k$ as well. Fix points $b_i\in U_i \cap U_{i-1}\cap \a_1$ and $a_i \in  U_i \cap U_{i-1}\cap \a_R$, so that the pairs $(a_i, b_i)$ are  endpoints of the same tunnel. Note that for every $i\in \ZZ_{k+1}$ an oriented quadrilateral loop $Q_i$ with vertices $[b_i,b_{i+1}, a_{i+1},a_i]$ (where each edge is an appropriate simple path contained in $\a_1, \a_R$ or in some tunnel) is contained in $U_i$ as suggested by Figure \ref{spaniergroupZ}. Define a
path $c_i$ between $x_0$ and $b_i$ to be the restriction of $\a_1$ to the appropriate interval. The loop
$$
(c_k ~ Q_k ~ c_k^{-1}) ~ (c_{k-1} ~ Q_{k-1} ~ c_{k-1}^{-1})\ldots (c_0 ~ Q_0 ~ c_0^{-1})
$$
of the Spanier group with respect to $\calU$
is homotopic to the non-trivial $\a_1$. Since $\a_1$ does not depend on a cover $\calU$, the based Spanier group of $Z'$ is non-trivial.

(5) Implied by (6) and Remark \ref{homHausRelation}.

(6) This is proved in \cite[Theorem 3.5]{CMRZZ}.

(7) We adopt the notation of part (1). Let $\a_1$ be a simple loop in the outer cylinder based at $x_0$.
Assume that $x_0$ is the endpoint of some tunnel, and that $\tau_R$ is the segment of this tunnel of length
$1\mathord{-}R$
starting at $x_0$. We will only use such $R$-values where the endpoint of $\tau_R$
is an intersection point with the surface portion.
We  prove that the conditions of Definition \ref{homHaus}(3) cannot be satisfied for
$w=\alpha_1$ and $v=\tau_R*\alpha_R*\tau_R^{-1}$. Assume the opposite. Then,
since $w\not\simeq1$ but $v\simeq1$ independently of $R$,
there is a finite covering $U_1, \ldots, U_k$ of $\a_1([0,1])$ as in Definition \ref{homHaus}(3). Without loss of generality we can decrease the sets $U_i$ to  path-connected sets, so that also the intersections $U_1\cap U_k$
and $U_i\cap U_{i+1}$ for $i=1,\ldots k-1$ are path-connected. Let $R<1$ be big enough so that $U_1, \ldots, U_k$ covers $\tau_R*\alpha_R*\tau_R^{-1} = v$.  We can reparametrize $v$ so that $v([t_{j-1},t_j])\subset U_j, \forall j=1,\ldots,k$. For every $j<k$ define $c_j$ to be a path in $U_j \cap U_{j+1}$ between $w(t_j)$ and $v(t_j)$, and for simplicity let $c_0=c_k$
be the constant path at $x_0$. Now for $i=1,\ldots,k$ let $w_i := c_{i-1}* v|_{[t_{i-1},t_i]}* c_i^{-1}$.
Obviously $w_1*\ldots* w_k\simeq v$, contradicting the conclusion of Definition \ref{homHaus}(3).
\hfill $\blacksquare$

\begin{defn}
\label{ShapeInj}
We call a space $X$ {\rm shape injective}, if the natural
homomorphism $\pi_1(X) \to \check{\pi}_1(X)$ is injective. For further information
see Section 3 of \cite{FZ}.
\end{defn}

\begin{prop}\label{HEprop}
The (based and unbased) Spanier groups
of the Hawaiian Earring are trivial.
\end{prop}

{\bf Proof.} We argue for the unbased Spanier-group, since
the based Spanier group is just a subgroup of it, anyhow.

This proposition can be deduced from the literature, namely from the
shape-injectivity
that the Hawaiian Earring enjoys as a
one-dimensional set (\cite[Theorem 1.1]{EK1} or \cite[Theorem 5.11]{CC})
or as planar set (\cite[Theorem 2]{FZ1}), together with the proposition \cite[4.8]{FZ}
which in this case implies that the Spanier group
is contained in a kernel of an injective map.
\hfill $\blacksquare$

\begin{theo}\label{hard}
The space $Y'$ is homotopically path-Hausdorff.
\end{theo}

{\bf Sketch of the proof.}
A space can only fail to fulfill the homotopic path-Hausdorff\-ness, if
arbitrarily close to some path $w$ there exist homotopic representatives
of the same path $v$ with $v \not\simeq w$. In such a case the space must
be wild in the neighbourhood of the trace of $w$.
Now $Y'$  is wild only at the central
axis, which  is a contractible part of the space. Thus $w$ will be in the trivial homotopy class,
and $v$ must be non-trivial. Therefore $v$ must leave the central axis through some
tunnel and return in such a way that it cannot be deformed onto the
central axis (e.g., by using other tunnels, or since it has
circled around the central axis before return). In any case $v$ must have
left the tunnel and also spent a segment on the surface.
 However, the place where $v$ leaves the tunnel and
continues via running over the waved surface, will be a characteristic of
the non-trivial homotopy class of $v$. Now, when in order to violate the
condition of Definition \ref{homHaus}(3), we have to construct another homotopic
representative of $v$ subordinated to a covering that does not contain
the places mentioned in the previous sentence, it will
belong to a different homotopy class. Thus our attempts to violate
Definition \ref{homHaus}(3) are bound to fail.\par
In a projected
forthcoming publication we plan to extend the techniques of
combinatorially describing homotopy classes of paths
in $H\!\!E   $
from \cite[2.3--2.10]{Z-Constr}, \cite[Section 1]{Z-Sp} to the
space $Y'$. Based on the appropriate combinatorial tool we plan to publish
a precise proof of this theorem, also. \par
The situation (cf.\ Proposition \ref{Z'}(7)) was different for the space $Z'$:
Here the wild part contained non-nullhomotopic curves;
thus it was possible that $w\not\simeq1$ and $v\simeq1$.
Indeed all paths $\tau_R*\alpha_R*\tau_R^{-1}$
that we constructed in the proof of Proposition \ref{Z'}(7)
were nullhomotopic.

\section{Proof of Theorem \ref{properties}: Algebraic characterization of semilocal
simple connectivity}

(1, ``$\Leftarrow$''):
If the Spanier group $\pi(\calU,x_0)$ vanishes for some covering $\calU$, then the products $u_1 v_1 u_1^{-1}$ of expression (\ref{UnbasedSpanierGroup}) are nullhomotopic loops. Such a product  is contractible if and only if $v_1$ is
contractible as both loops are freely homotopic. Furthermore, the contractibility of a loop  never depends on whether
the endpoint is kept fixed.
Thus the elements of the covering $\calU$ suffice to prove that $X$ is unbased semilocally simply connected as every loop $v_1$ that is contained in some element of $\calU$ is contractible.

(1, ``$\Rightarrow$''):
For every point of $X$ choose
a neighbourhood which satisfies the condition of Definition \ref{UnbasedSLSC} and form
the covering $\calU$  of these neighbourhoods. Every product in the form of expression (\ref{UnbasedSpanierGroup})
can be contracted by first contracting the $v_i$-loops, and then contracting $u_i ~ u_i^{-1}$. Hence the Spanier group is trivial.

(2, ``$\Leftarrow$''): {If for some open covering  $
\calV = \{(U_j,x_j) \mid j \in J\}$ by pointed sets $(U_j,x_j)$ the Spanier group $
\pi^*(\calV,x_0)$ vanishes}, then also the products $u_1 v_1 u_1^{-1}$ of expression (\ref{BasedSpanierGroup}) are nullhomotopic paths. Based on these
substitutions, the same line of arguments as in (1, ``$\Leftarrow$")  can be used. However, while in expression (\ref{UnbasedSpanierGroup})
$v_i$ could have been any loop in a neighbourhood $U_j$, the definition of $\pi^*$ requires to consider only such $v_i$
that are based at $x_j$. If $U_j$ should contain loops whose trace has no path-connection to $x_j$
 inside $U_j$, our assumptions do not suffice to conclude that such loops will also be contractible. Thus in this case we can only conclude that Definition \ref{BasedSLSC}, but not that Definition \ref{UnbasedSLSC} will be fulfilled.

(2, ``$\Rightarrow$''):
For every point $x\in X$ choose a neighbourhood $U_x$ so that $\pi_1(U_x,x)\to \pi_1(X,x)$ is zero. {Form an open covering $\calU$  of pointed sets $\{(U_x,x)\mid x\in X\}$ and consider its based Spanier group}. Every product in the form of expression (\ref{BasedSpanierGroup})
can be contracted by first contracting the loops $v_i$, and then contracting $u_i ~ u_i^{-1}$. Hence the based Spanier group is trivial.

(3): It suffices to observe that the triviality of all loops implies triviality of all based loops within any arbitrary neighbourhood.

(4): Given a point $x\in X$ choose a neighbourhood $U$ of $x$ that satisfies the condition of Definition \ref{BasedSLSC}. Then any path-connected neighbourhood $V\subset U$ of $x$ satisfies the condition of Definition \ref{UnbasedSLSC}.

(5): See Proposition \ref{Y}. \hfill $\blacksquare$
\bigskip

\section{Proof of Theorem \ref{covering} : Algebraic criteria for homotopic Hausdorffness}

{\bf Convention:} Within this section the overline (``$\overline{(~)}$")
denotes the reversion of the orientation of a path.

{\bf Proof: }We split the statement of Theorem \ref{covering} into two implications:

$\invlim\pi^*({\mathcal U})$ is trivial.
$\buildrel(1)\over\Longrightarrow$
  homotopic path-Hausdorffness
$\buildrel(2)\over\Longrightarrow$
             existence of a generalized universal covering space.

(1): Suppose $X$ is not homotopically path-Hausdorff. Using the notation of Definition \ref{homHaus}(3) there exist paths $w,v \colon [0,1]\to X$ with $w(0)\mathord{=}v(0)\mathord{=}P\Komma w(1)\mathord{=}v(1)\mathord{=}Q$ and a non-trivial homotopy class  $\a\in \pi_1(X,P)$, $\alpha :=[ w*\overline{v}]$ for which the conditions of Definition
\ref{homHaus}(3) do not hold.  We claim that $\a$ is contained in $\invlim\pi^*({\mathcal U})$. The proof will resemble that of Proposition \ref{Z'}(1){.

Let $\mathcal U$ be an open cover of $X$ by pointed sets
(cf.\ Definition \ref{nbd-pair}). Choose a cover $U_1, U_2, U_3,\ldots, U_k$  of $w([0,1])$ by
 open sets from $\mathcal U$  so that there exists a  partition $0=t_0 < t_1 < t_2\ldots < t_k=1$ for which  $U_j$ covers
$w([t_{j-1},t_j])$ and the according base point for each of the $U_j$ lies on
the segment $w([t_{j-1},t_j])$. Since Definition \ref{homHaus}(3) is not satisfied there exist paths
$w_j$ such that $w_j$  connects $w(t_{j-1}) $ with $w(t_j)$ inside $U_j$, and so that the concatenation
$v':= w_1*\ldots*w_k$ is homotopic to $v$. Assume that $v'$ is parametrized so that $v(t_j) = w(t_j)$
for $j=0,\ldots,k$.
Note that the concatenation
$$
\big({v|_{[t_0,t_0]}}*w|_{[t_0,t_1]} * \overline{v|_{[t_0,t_1]}}\big) * \big({v|_{[t_0,t_1]}} * w|_{[t_1,t_2]} * \overline{v|_{[t_0,t_2]}} \big) * \ldots
$$
$$
\ldots* \big({v|_{[t_{0},t_{k-1}]}} * w|_{[t_{k-1},t_k]} * \overline{v|_{[t_0,t_k]}} \big)
$$
is  homotopic to $\a$ and contained in $\pi^*({\mathcal U},P)$ as each of the factors ${v|_{[t_{0},t_{i-1}]}} * w|_{[t_{i-1},t_i]} * \overline{v|_{[t_0,t_i]}}={v|_{[t_{0},t_{i-1}]}} * w|_{[t_{i-1},t_i]} * \overline{v|_{[t_{i-1},t_i]}}*\overline{v|_{[t_{0},t_{i-1}]}}$ is a conjugate of the loop $w|_{[t_{i-1},t_i]} * \overline{v|_{[t_{i-1},t_i]}}$, that is
contained in  $U_j\in \mathcal U$ and has its base point $w(t_{i-1})$ inside $U_j$ path-connected to the
base point of $U_j$.}
\medskip

(2): We give a proof by contradiction. Recall that the standard
existence proof of covering spaces (e.g.\ \cite[p.\ 393]{Mu}) is based on
interpreting the Universal Path Space (as it was called by \cite{BS})
as a covering space. Also, the generalized universal covering spaces in the
sense of \cite{FZ} are constructed on the basis of considering
the Universal Path Space, with an adaptation of the definition
of the topology to the situation of the absence of
semilocal simple connectivity. Therefore the points in the covering space of $(X,x_0)$
are represented by (homotopy classes of) paths in the base space and the topology is induced by the sets
$$
\UU(\g,U)=\{\g* \d; \d\colon [0,1]\to U, \d(0)=\g(1)\}
$$
where $\g\colon [0,1]\to X$ is a path  originating at $x_0$ and $U\subset X$ is an open neighbourhood of $\g(1)$.

The assumption for the desired proof by contradiction is the following:
There does not exist a generalized universal covering space. By \cite[2.14]{FZ}
this means,\ that there exists a path $w\colon [0,1]\to X$ which allows two different
lifts to the covering space with the same start-point. We can assume that $t\mapsto w_t:=w|_{[0,t]}$ and $t\mapsto v_t$ are two different lifts of path $w$ with $v_0=w_0$ being a constant path at $w(0)$ and $v := v_1 \not\simeq w_1=w$. We claim that such a situation does not comply with the conditions of Definition \ref{homHaus}(3).

\begin{figure}
\includegraphics[scale=1.176470]{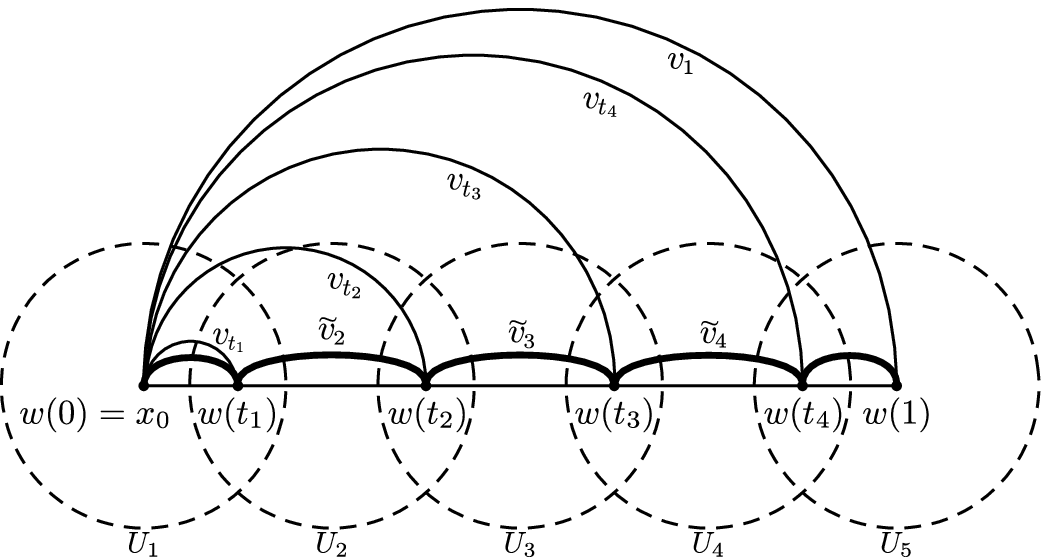}
\caption{Two different lifts of the path $w$.\newline
         The drawing illustrates the proof-construction of step (2) in the
         special situation where we have $k\mathord{=}5$. In addition, this figure is drawn
         in a way that for each $j$ we just get $l=1$, i.e.\
         just the case of a trivial
         $s$-partition is pictured. The line in bold corresponds to the path $
         \tilde v$.
}\label{NonUniqueLift}
\end{figure}

Choose any covering $U_1,\ldots,U_k $ of $w([0,1])$
such that for a suitable partition $0=t_0 < t_1 < t_2 < t_3 <\ldots< t_k=1\Komma
U_j$ covers $w([t_{j-1},t_j])$, as suggested by Figure \ref{NonUniqueLift}. Let us focus on the interval $[t_{i-1},t_{i}]$ for some fixed $i$. By the continuity of the lift $v$ there exists for every $t\in [t_{i-1},t_{i}]$ a neighbourhood $V_t$ of $t$ so that $v_s\in \UU(v_t,U_i), \forall s\in V_t$. Hence we can find a partition $t_{i-1}=s_0 < s_1 < \ldots < s_l=t_i$ so that $v_s\in \UU(v_{s_j},U_i), \forall s\in [s_{j},s_{j+1}]$, i.e.\ $v_{s_j}*\d_j \simeq v_{s_{j+1}}$ for some path $\d_j$ in $U_i$. This condition implies the existence of a path $\tilde v_i\colon [t_{i-1},t_{i}]\to U_i$ between $v_{t_{i-1}}(1)$ and $v_{t_i}(1)$, so that $v_{t_{i-1}}* \tilde v_i \simeq v_{t_i}$. The paths $\tilde v_i$, defined on $[t_{i-1},t_{i}]$, induce a path $\tilde v\colon [0,1] \to X$.

Note that $$\tilde v := \tilde v_1 * \tilde v_2 * \ldots * \tilde v_k\simeq (\overline{v_0}*v_0)*\tilde v_1 * (\overline{v_{t_1}} * v_{t_1}) *
\tilde v_2 *
\ldots *(\overline{v_{t_{k-1}}} * v_{t_{k-1}}) *\tilde v_k\simeq$$$$\simeq
({v_0} *
\tilde v_1 )* \overline{v_{t_1}} * (v_{t_1} *
\tilde v_2) *\overline{v_{t_2}} * (v_{t_2}*\tilde v_3)* \ldots *\overline{v_{t_{k-1}}} * (v_{t_{k-1}} *\tilde v_k) \simeq
$$$$\simeq
v_{t_1} * \overline{v_{t_1}} * v_{t_2} *\overline{v_{t_2}} * v_{t_3}
*\ldots * \overline{v_{t_{k-1}}} *v_{t_{k}} \simeq
v_{t_k} =
 v_1 = v$$ hence
with letting $w_i := \tilde v_i$ we obtain that Definition \ref{homHaus}(3) is not fulfilled.\hfill $\blacksquare$
\bigskip

\section{Overview of the implications}

\begin{figure}[htbp]

\epsfig{file=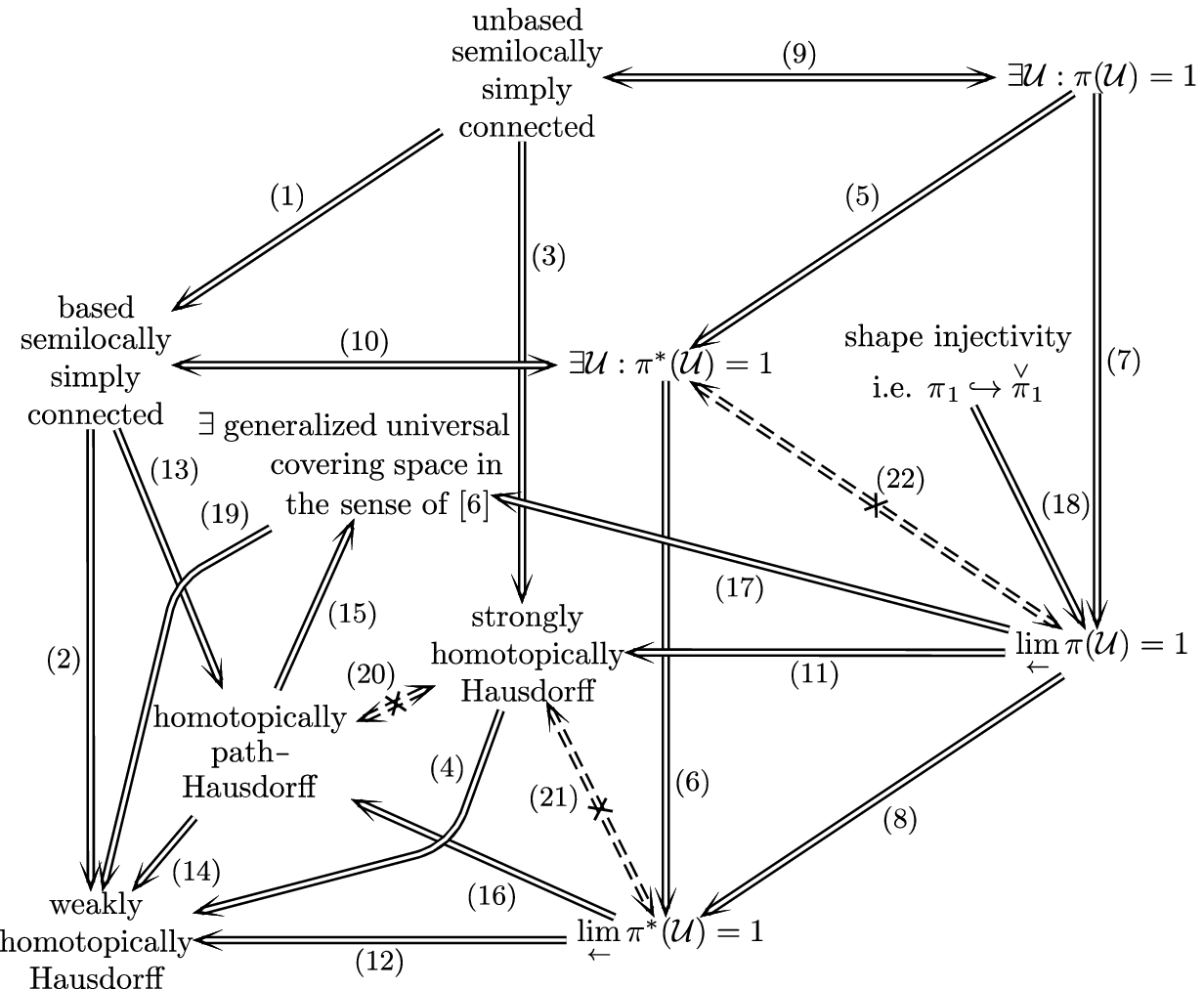,height=4in}

\caption{{ Diagram}}
\label{diagram}
\end{figure}

The diagram of this section (see Fig. \ref{diagram}) gathers together a number of implications
of properties of a space
that occurred
in, or are closely related to, the ones from our Theorems \ref{properties} and \ref{covering}.
No  assumption other than path-connectedness (cf.\ Remark \ref{remark}) is made here.
According to the enumeration of the implications in the diagram, for each arrow a reference or a sketch of the proof
is given. The label (1, ``$\Longrightarrow$'') means, that an argument is to be given,
why this implication is true, while (1, ``\drauf{$\Longleftarrow$}{\Schrstr}'') means,
that an argument is to be given, why the converse of this implication is in general
not true.

(1, ``$\Longrightarrow$''):  Follows from the definition of $\pi_1$.

(1, ``\drauf{$\Longleftarrow$}{\Schrstr}''):  The counterexample is $Y$ (see Proposition \ref{Y}).

(2, ``$ \Longrightarrow$''): This is a passage to an obviously weaker property.

(2, ``\drauf{$\Longleftarrow$}{\Schrstr}''):  The Hawaiian Earring is an appropriate example. It contains arbitrarily small essential loops in the neighbourhood of the accumulation point. Therefore it is not semilocally simply connected.  On the other hand,  one-dimensional spaces are weakly homotopically Hausdorff by \cite[Corollary 5.4(2)]{CC}.

(3, ``$\Longrightarrow$''): The assumption means that every point has a neighbourhood such
   that all loops in this neighbourhood are contractible. Such  neighbourhoods suffice to prove that the space is strongly homotopically Hausdorff.

(3, ``\drauf{$\Longleftarrow$}{\Schrstr}''):  Similarly as for (2, ``\drauf{$\Longleftarrow$}{\Schrstr}'') the Hawaiian Earring is an appropriate example. It is easy to see that it is not unbased semilocally simply connected at the accumulation point due to arbitrarily small essential loops. It is slightly harder to see that it is strongly homotopically Hausdorff.

Let us denote the accumulation point of the Hawaiian Earring $H\!\!E$  by $x_0$ and let us enumerate the loops of $H\!\!E$ by $l_i, i=1,2,\ldots$, i.e.\ $\cap_i l_i =\{x_0\}$
(cf.\ Figure~\ref{HEfig}).
It is enough to consider the condition for strong homotopic Hausdorffness at $x_0$ as all other points have contractible neighbourhoods. Let $\{V_i\}_{i=1,2,\ldots}$ denote a basis of open neighbourhoods of $x_0$ so that $l_i\subset V_j$ iff $i \geq j$. Suppose there exists a loop $\a=\a_1$ which is freely homotopic to some loop $\a_i$ in $V_i$ for every $i>1$. Note that there is a strong deformation retraction
$$
V_i \to W_i:=\bigcup_{j\geq i}l_j.
$$
Therefore we can assume that every $\a_i$ is a loop in $W_i$. For every $i$ the loops $\a_i$ and $\a_{i+1}$ are homotopic. By Lemma 4.3(1) of \cite{CMRZZ} the homotopy $h_i$ between them can be chosen within $W_i$. We define a map $f\colon B^2(o,1) \to H\!\!E$ on a closed unit disc
with midpoint $o$ by assigning
$f|_{S^1(o,1/i)}=\a_i$  and by using the  homotopies $h_i$ for
defining
$$
f|_{\overline{{B^2(o,1/i)-B^2(o,1/(i+1))}}}:=h_i
$$
on the closed annuli
between the concentric circles $S^1(o,1/i)$.
The map $f$ is obviously continuous on $B^2(o,1)-\{o\}$ as it is obtained by a locally finite gluing of maps that agree on the intersections of annular domains (i.e.\ on concentric circles). It is also continuous at $o$ because for every $i$ the  preimage of $V_i$ contains  $B^2(o,1/i)$. The map $f$ provides a nullhomotopy for $\a=\a_i$.

(4, ``$\Longrightarrow$''): This is a passage to a weaker requirement.

(4, ``\drauf{$\Longleftarrow$}{\Schrstr}''): The space $Y$ serves as an example by Proposition \ref{Y}.

(5, ``$\Longrightarrow$''):  This is a weakening of the conditions, since
    $\pi^*$-groups
   contain fewer elements than  $\pi$-groups.

(5, ``\drauf{$\Longleftarrow$}{\Schrstr}''):  The space $Y$ serves as an example. To prove it we use Proposition \ref{Y}. The non-triviality of the Spanier group of $Y$ and (7) assure that there is no cover with non-trivial Spanier group.\ Yet $Y$ is semilocally simply connected.

(6, ``$\Longrightarrow$''):  By Remark \ref{remark0} the condition ``$\pi^*(\mathcal U)=1$"
    for some cover $\mathcal U$    guarantees the triviality of the inverse limit.

(6, ``\drauf{$\Longleftarrow$}{\Schrstr}''):  Again the Hawaiian Earring gives the corresponding example. Every covering of the accumulation point contains small essential loops entirely and thus no $\pi^*(\mathcal U)$-group can be trivial.
In Proposition \ref{HEprop} it was shown that it has trivial Spanier group.

(7, ``$\Longrightarrow$''):  The argument  here is analogous
    to that of (6, ``$\Longrightarrow$'').

(7, ``\drauf{$\Longleftarrow$}{\Schrstr}''):  The same example as for (6, ``\drauf{$\Longleftarrow$}{\Schrstr}'') suffices.
    For locally path-connected
    spaces there is no difference between the
    based and unbased statements.

(8, ``$\Longrightarrow$''): Follows from the same argument as (5, ``$\Longrightarrow$'').

(8, ``\drauf{$\Longleftarrow$}{\Schrstr}''): The space $Y$ can be used for this purpose by Proposition \ref{Y}.

(9, ``$\Longleftrightarrow$''):  This is the statement of Theorem \ref{properties}(1).

(10, ``$\Longleftrightarrow$''):  This is the statement of Theorem \ref{properties}(2).

(11, ``$\Longleftarrow$''):   If a space is not strongly homotopically
Hausdorff, then  there exists a point  $P$ and a non-trivial free homotopy class $\alpha$ which
can represented in any neighbourhood of $P$. This class will appear in every
$\pi(\mathcal U)$ (and thus in their inverse limit), since each $\mathcal U$
has to contain a neighbourhood of $P$.

(11, ``\drauf{$\Longrightarrow$}{\Schrstr}''):  The corresponding example is $Z'$. By Proposition \ref{Z'} it is strongly homotopically Hausdorff but has a non-trivial Spanier group.\

(12, ``$\Longleftarrow$''):  The same argument as for (11) applies in terms of based homotopies.

(12, ``\drauf{$\Longrightarrow$}{\Schrstr}''):  The same counterexample as for (11) applies as the space $Z'$ is locally path-connected, thus the base points do not matter.

(13, ``$\Longrightarrow$''):   The implication  follows from  statements (10), (6) and (16).

(13, ``\drauf{$\Longleftarrow$}{\Schrstr}''):
The Hawaiian Earring can  be used  as an example. It is obviously  not semilocally simply connected. On the other hand, the proof of (6, ``\drauf{$\Longleftarrow$}{\Schrstr}'') establishes the triviality of the based Spanier group of $H\!\!E$, hence it is homotopically path-Hausdorff by (16).

(14, ``$\Longrightarrow$''):  This  has been observed in Remark
\ref{homHausRelation} already.

(14, ``\drauf{$\Longleftarrow$}{\Schrstr}''): The space $Z'$ can be used for this purpose by Proposition \ref{Z'}.

(15, ``$\Longrightarrow$''):
This is implication (2) from the proof of Theorem \ref{covering} in Section 5.

(16, ``$\Longrightarrow$''):  This is implication (1) from the proof of Theorem \ref{covering} in Section 5.

(16, ``\drauf{$\Longleftarrow$}{\Schrstr}''):   The corresponding example  is $Y'$ by Proposition \ref{Y'}.

(17, ``$\Longrightarrow$'') and (18, ``$\Longrightarrow$''):
These implications have been proved in  \cite[4.7]{FZ}  and \cite[4.8]{FZ}  respectively. Implication (17) also follows from implications (8),(16) and (15).

(17, ``\drauf{$\Longleftarrow$}{\Schrstr}''):
 The corresponding example  is $Y'$:
By Proposition \ref{Y'}(1) it has non-trivial Spanier group, but according
to (15) and Proposition \ref{Y'}(7)  it has a generalized covering space.

(18, ``\drauf{$\Longleftarrow$}{\Schrstr}''):
The corresponding example  is $Z$. The fundamental group of $Z$ is   $\ZZ$ generated by the outer cylinder.
This follows, since by \cite[Lemma 4.3(1)]{CMRZZ} any nullhomotopy of a loop that is not passing through the arc $C$
need not to pass through $C$, either, but $Z - C$ consists of two different path-components:
the outer cylinder, and the contractible surface portion.
The fundamental  group of $Z$ vanishes when we pass to the shape group (i.e.\ the outer cylinder is homotopically trivial in every neighbourhood of $Z$ in $\RR^3$), hence the space is not shape injective. On the other hand the Spanier group of $Z$ is trivial by Proposition \ref{Z}(1).

(19, ``{$\Longleftarrow$}''): By contradiction: The presence of small loops in the sense of \cite{ZV} is equivalent to the absence of weak homotopic Hausdorffness. Every small loop induces a non-trivial lift by (cf.\ \cite[Lemma 16]{ZV} or \cite[Lemmas 2.10--2.11]{FZ}),
hence there is no generalized universal covering space in terms of \cite{FZ}.

(20, ``\drauf{$\Longleftrightarrow$}{\Schrstr}''):  The corresponding examples  are $Y$  and $Z'$ by Propositions \ref{Y} and \ref{Z'}.

(21, ``\drauf{$\Longleftrightarrow$}{\Schrstr}''):
The corresponding examples  are $Y$  and $Z'$ by Propositions \ref{Y} and \ref{Z'}.

(22, ``\drauf{$\Longrightarrow$}{\Schrstr}''):
The corresponding example is  $Y$ by Proposition \ref{Y}(1) and
the proof of Proposition \ref{Y}(2).

(22, ``\drauf{$\Longleftarrow$}{\Schrstr}''):
The corresponding example is  the Hawaiian Earring.
The argument is the same as that of   ``(6 \drauf{$\Longleftarrow$}{\Schrstr}'',
since there is no difference between based and unbased statements for locally path-connected
    spaces.\hfill $\blacksquare$

\section{Acknowledgements}

This research was supported by the Polish--Slovenian  grants BI-PL 2008-2009-010
and 2010-2011-001,
the second and the third authors were supported by
the ARRS program P1-0292-0101 and
project J1-2057-0101, and the fourth author was supported by
KBN grant N200100831/0524.
The authors
wish to thank
the referee for valuable
comments and suggestions.
\par

\end{document}